\documentclass[a4paper,12pt]{article}
\usepackage{amsmath}
\usepackage{amssymb}
\usepackage{tabularx}
\usepackage{enumerate}
\usepackage[numbers,sort&compress]{natbib}

\newtheorem{thm}{Theorem}[section]
\newtheorem{lem}[thm]{Lemma}

\newtheorem{rmk}{Remark}[section]
\newtheorem{pppp}{Proof}

\newcommand{\qed}{\hspace{1em}\mbox{\raisebox{0.65ex}{\fbox{}}}}

\numberwithin{equation}{section}
\newcommand{\ol}{\overline}

\newcommand{\be}{\begin{equation}}
\newcommand{\ee}{\end{equation}}
\newcommand\bes{\begin{eqnarray}} \newcommand\ees{\end{eqnarray}}
\newcommand{\bess}{\begin{eqnarray*}}
\newcommand{\eess}{\end{eqnarray*}}

\newcommand{\R}{\mathbb{R}}
\newcommand{\bpf}{{\bf Proof:\ \ }}
\newcommand{\epf}{\mbox{}\hfill $\Box$}

\begin{document}
\setlength{\baselineskip}{2\baselineskip}
\thispagestyle{empty}

\title{An SIR epidemic model with free boundary\thanks{The work is supported by BSR Program of NRF/MEST (Grant N0.
2010-0025700), the PRC grant NSFC 11071209, and also by the Ph.D. Programs Foundation of Ministry of
Education of China (No. 20113250110005).}}
\date{\empty}

\author{Kwang Ik Kim$^{a}$, Zhigui Lin$^b$\thanks{Corresponding author. Tel. +86 514 87975401, Fax +86 514 87975423.
 Email: zglin68@hotmail.com (Z. Lin).} and Qunying Zhang$^b$\\
{\small $^a$Department of Mathematics, Pohang University of Science and Technology, } \\
{\small Pohang, 790-784, Republic of Korea. }\\
{\small $^b$School of Mathematical Science, Yangzhou University, Yangzhou 225002, China.}\\
}

 \maketitle

\begin{quote}
\noindent
{\bf Abstract.} { 
\small An SIR epidemic model with free boundary is investigated. This model
describes the transmission of diseases. The behavior
of positive solutions to a reaction-diffusion system in a radially symmetric domain is investigated.
The existence and uniqueness of the global solution are given by the contraction mapping theorem.
Sufficient conditions for the disease vanishing or spreading are given. Our
result shows that the disease will not spread to the whole area if the basic reproduction
number $R_{0}<1$ or the initial infected radius $h_0$
 is sufficiently small even that  $R_{0}>1$. Moreover, we prove that the disease will spread to the whole area
if $R_{0}>1$ and the initial infected radius $h_0$ is suitably large.}

\noindent {\it MSC:} primary: 35B35; secondary: 35K60

\medskip
\noindent {\it Keywords: } Reaction-diffusion systems; SIR model;
Free boundary; Dynamics

\end{quote}

\section{Introduction}

Recently epidemic model has been received a
great attention in mathematical ecology. To describe the development of an infectious disease,
compartmental models have been given to separate a population into various
classes based on the stages of infection \cite{AM3}. The classical SIR model
is described by partitioning the population into
susceptible, infectious and recovered individuals, denoted by $S, I$ and $R$, respectively. Assume that the disease
incubation period is negligible so that each susceptible
individual becomes infectious and later recovers with a permanently
or temporarily acquired immunity, then the SIR model is governed by the following system of differential
equations:
\bes \left\{
\begin{array}{lll}
\dot{S}(t)&=&-\beta S(t)I(t)-\mu _1S(t)+b,\\[1mm]
\dot{I}(t)&=&\beta S(t)I(t)-\mu _2I(t)-\alpha I(t),\\[1mm]
\dot{R}(t)&=&\alpha I(t)-\mu _3R(t),
\end{array} \right.
\label{a1} \ees
where the total population size has been normalized to one and the influx of the
susceptible comes from a constant recruitment rate $b$. The death rate for the $S, I$ and $R$
class is, respectively, given by $\mu_1, \mu_2$ and $\mu_3$.
Biologically, it is natural to assume that $\mu_1<\min \{\mu_2, \mu_3\}$.
 The standard incidence of disease is denoted by $\beta SI$, where
$\beta$ is the constant effective contact rate, which is the average
number of contacts of the infectious per unit time.
The recovery rate of the infectious is denoted by $\alpha$ such that $1/\alpha$ is the mean
time of infection.

In \cite{6}, the threshold behavior was given. The authors showed that the basic reproduction
number $R_0$ ($=\frac{b \beta}{\mu_1(\mu_2+\alpha)}$) determines
whether the disease dies out ($R_0<1$) or remains endemic
($R_0>1$). In \cite{4}, a complete analysis of the global
dynamics of an ordinary differential equation model with multiple infectious stages was presented, showing the same
threshold behavior. For other works on various types of SIR
epidemic model, interested readers may refer to \cite{AM1, BHM, Gr, Ji, LM, MLL, ZH} and the references therein.

There are other compartmental combinations for modelling some other diseases. For example,
 the SI model describes a disease, such as herpes or HIV, with two stages, where
individuals are infectious for life and never removed. The SIS model describes the case
when individuals recover from the disease but there is no immunity, and they return
to the susceptible class. Examples for this SIS model include sexually transmitted diseases, plague and
meningitis. Unlike SIR models, SEI models  \cite{GMH, LZ1} assume that a susceptible
individual first goes through a latent (exposed) period before
becoming infectious. An example of this model is the transmission
of SARS \cite{XJ}, which is one of the serious diseases that human
beings face at present.

When the distribution of the distinct classes is in different spatial
locations, the diffusion terms should be taken into consideration
and thus an extended version of the above SIR system (\ref{a1}) can be described as the following:
\begin{eqnarray}
\left\{
\begin{array}{lll}
S_{t}-d_1\Delta S=-\beta S(t)I(t)-\mu _1S(t)+b,\; &
x\in\Omega,\; t>0,  \\[1mm]
  I_{t}-d_2\Delta
I=\beta S(t)I(t)-\mu _2I(t)-\alpha I(t),\; &
x\in\Omega,\; t>0, \\[1mm]
R_{t}-d_3 \Delta R=\alpha I(t)-\mu _3R(t),\; &
x\in\Omega,\; t>0, \\[1mm]
\partial_\eta S=\partial_\eta
I=\partial_\eta R=0\; &x\in\partial\Omega,\; t>0,\\[1mm]
S(x,0)=S_{0}(x),\, I(x,0)=I_{0}(x),\,R(x,0)=R_{0}(x),\, &x\in \ol \Omega ,
\end{array} \right.
\label{a2}
\end{eqnarray}
where $\Omega$ is a fixed and bounded domain in $\R^n$ with smooth
boundary $\partial \Omega$, and $\eta$ is the outward unit normal
vector on the boundary. Here the homogeneous Neumann boundary
condition implies that the above system is self-contained and
there is no emigration across the boundary. The positive constants
$d_i(i=1,2,3)$ are the diffusion coefficients.

It must be pointed out that the solution of system (\ref{a2}) is always positive for any time $t>0$  no matter what the nonnegative
nontrivial initial date is. It means that the disease spreads to the whole area immediately even when
the infectious is confined to a small part of the area in the beginning. It doesn't match the
observed fact that disease always spreads gradually. Recently the free boundary has been introduced in many areas, especially
 the well-known Stefan condition has been used to describe the spreading process. For example, it was used
in describing the melting of ice in contact with water \cite{R}, in
the modeling of oxygen in the muscle \cite{C}, and in the dynamics
of population \cite{HMS, KLL, LIN, MYY}. There is a vast literature
on the Stefan problem, and some recent and theoretically advanced results can be found in \cite{CS}.

Motivated by the statements mentioned above, we are attempting to consider a SIR epidemic model
with a free boundary, which describes the spreading frontier of the disease.  For simplicity, we
assume the environment is radially symmetric. We will investigate the behavior of the positive solution
$(S(r, t), I(r, t), R(r, t); h(t))$ with $r=|x|$ and $x \in \R^n$ in the following problem:
\begin{eqnarray}
\left\{
\begin{array}{lll}
S_{t}-d_1\Delta S=b-\beta S(r, t)I(r, t)-\mu _1S(r, t),\; &
r>0,\; t>0,  \\[1mm]
  I_{t}-d_2\Delta
I=\beta S(r, t)I(r, t)-\mu _2I(r, t)-\alpha I(r, t),\; &
0<r<h(t),\; t>0, \\[1mm]
R_{t}-d_3 \Delta R=\alpha I(r, t)-\mu _3R(r, t),\; &
0<r<h(t),\; t>0, \\[1mm]
S_r(0, t)=I_r(0,t)=R_r(0,t)=0,&t>0,\\[1mm]
I(r,t)=R(r,t)=0,&r\geq h(t),\, t>0,\\[1mm]
h'(t)=-\mu I_{r}(h(t),t),\; h(0)=h_0>0,  & t>0, \\[1mm]
S(r,0)=S_{0}(r),\; I(r,0)=I_{0}(r),\; R(r, 0)=R_{0}(r), & r\geq 0,
\end{array} \right.
\label{a3}
\end{eqnarray}
where $\triangle w=w_{rr}+\frac{n-1}{r}w_r$, $r=h(t)$ is the moving
boundary to be determined,  $h_0,\, d_i$ and $\mu$ are positive constants. The initial functions
$S_0, I_0$ and $R_0$ are nonnegative and satisfy
\begin{eqnarray}
\left\{
\begin{array}{ll}
S_0\in C^{2}([0, +\infty)),\, \, I_0,\, R_0\in C^2([0, h_0]), \\[1mm]
I_0(r)=R_0(r)=0,\, r\in [h_0, +\infty) \ \textrm{and} \ I_0(r)>0,\, r\in [0, h_0).
\end{array} \right.
\label{Ae}
\end{eqnarray}
Ecologically, this model means that beyond the free boundary $r=h(t)$, there is only susceptible, no infectious or recovered individuals.
The equation governing the free boundary,
$h'(t)=-\mu I_{r}(h(t),t)$, is a special case of the well-known
Stefan condition, which has been established in \cite{LIN} for the diffusive populations.

The remainder of this paper is organized as follows.
In the next section, we first apply a contraction mapping theorem to
prove the global existence and uniqueness of the solution to the problem
(\ref{a3}).  Then we make use of the Hopf Lemma to give the
monotonicity of the free boundary. Section 3 is devoted to prove that the disease will vanish if the basic reproduction
number $R_{0}<1$.
In Section 4, we discuss the case $R_{0}>1$. Our results show that for the case $R_{0}>1$, the disease will spread to the whole area if $h_0$ is suitably large;
while the disease will vanish if $h_0$ is sufficiently small. Our
arguments are based on the comparison principle and the construction
of appropriate supper solution of (\ref{a3}). Finally, we give a brief discussion in Section 5.

\section{Existence and uniqueness}

In this section, we first prove the following local existence and
uniqueness result by the contraction mapping theorem. We then use
suitable estimates to show that the solution is defined for all
$t>0$.

\begin{thm} For any given $(S_0, I_0, R_0)$ satisfying \eqref{Ae} and any $\gamma \in (0, 1)$, there is a $T>0$ such that  problem \eqref{a3}
admits a unique bounded solution
$$(S, I, R; h)\in C^{1+\gamma, (1+\gamma)/2}(D^\infty_{T})\times[C^{ 1+\gamma,(1+\gamma)/2}(D_{T})]^2\times C^{1+\gamma/2}([0,T]);$$
moreover,
\begin{eqnarray}
\|S\|_{C^{1+\gamma,(1+\gamma)/2}({D}^\infty_{T})}+\|I\|_{C^{1+\gamma,
(1+\gamma)/2}({D}_{T})}+\|R\|_{C^{1+\gamma,
(1+\gamma)/2}({D}_{T})}+||h\|_{C^{1+\gamma/2}([0,T])}\leq C,\label{b12}
\end{eqnarray}
where $D^\infty_{T}=\{(r, t)\in \R^2: r\in [0, +\infty), t\in [0,T]\}$
and $D_{T}=\{(r, t)\in \R^2: r\in [0, h(t)], t\in [0,T]\}$. Here $C$
and $T$ only depend on $h_0, \gamma, \|S_0\|_{C^{2}([0, \infty))}, \|I_0\|_{C^{2}([0, h_0])}$ and $\|R_0\|_{C^{2}([0, h_0])}$.
\end{thm}
\bpf  We first straighten the free boundary as in \cite{CF}. Let
$\xi(s)$ be a function in $C^3[0, \infty)$ satisfying
$$\xi(s)=1\ \textrm{if} \,\, |s-h_0|<\frac{h_0}8,\; \xi(s)=0\
 \textrm{if} \,\, |s-h_0|>\frac {h_0}2,\ |\xi
'(s)|<\frac 5{h_0} \mbox{ for all } s.$$ Consider the transformation
$$(y, t)\rightarrow (x, t), \textrm{where}\,\, x=y+\xi (|y|)(h(t)-h_0y/|y|),
\quad y\in R^n,$$
which leads to the transformation
$$(s, t)\rightarrow (r, t)\ \textrm{with}\,\, r=s+\xi (s)(h(t)-h_0),
\quad 0\leq s<\infty.$$
As long as $$|h(t)-h_0|\leq \frac {h_0}8,$$
the above transformation $x \to y$ is a diffeomorphism from $R^n$
onto $R^n$ and the transformation $s\to r$ is also a diffeomorphism from $[0, +\infty)$
onto $[0, +\infty)$. Moreover, it changes the free boundary $r=h(t)$
to the line $s=h_0$. Now, direct calculations show that
\begin{eqnarray*}
\displaystyle\frac {\partial s}{\partial r}=
\frac1{1+\xi'(s)(h(t)-h_0)}&:=&\sqrt{A(h(t), s)},\\[1mm]
\displaystyle \frac {\partial^2 s}{\partial r^2}=-\frac
{\xi''(s)(h(t)-h_0)}{[1+\xi'(s)(h(t)-h_0)]^3}&:=&B(h(t),s),\\[1mm]
\displaystyle -\frac 1{h'(t)}\frac {\partial s}{\partial t}=\frac
{\xi(s)}{1+\xi'(s)(h(t)-h_0)}&:=&C(h(t), s),\\[1mm]
\displaystyle \frac {(n-1)\sqrt{A}}{s+\xi(s)(h(t)-h_0)}&:=&D(h(t), s).
\end{eqnarray*}

Now, if we set
$$S(r, t)=S( s+\xi (s)(h(t)-h_0),t):=u(s, t),$$
$$I(r, t)=I( s+\xi (s)(h(t)-h_0),t):=v(s, t),$$
$$R(r, t)=R( s+\xi (s)(h(t)-h_0),t):=w(s, t),$$
then the free boundary problem \eqref{a3} becomes
\begin{eqnarray}
\left\{
\begin{array}{lll}
u_{t}-Ad_1\Delta_s u-(Bd_1+h'C+Dd_1)u_s=b-\beta uv-\mu _1u,\; &s>0,\; t>0, \\[1mm]
v_{t}-Ad_2\Delta_s v-(Bd_2+h'C+Dd_2)v_s=\beta uv-\mu _2v-\alpha v,\; &0<s<h_0,\; t>0,\\[1mm]
w_{t}-Ad_3\Delta_s w-(Bd_3+h'C+Dd_3)w_s=\alpha v-\mu _3w,\; &0<s<h_0,\; t>0, \\[1mm]
u_s(0, t)=v_s(0, t)=w_s(0,t)=0, &t>0,\\[1mm]
v(s,t)=w(s,t)=0,&s\geq h_0,\, t>0,\\[1mm]
h'(t)=-\mu v_s(h_0,t),\,h(0)=h_0,&t>0, \\[1mm]
u(s,0)=u_{0}(s),\,v(s,0)=v_0(s),\,w(s,0)=w_0(s),&s\leq 0,
\end{array}\right.
\label{Bb}
\end{eqnarray}
where $A=A(h(t),s)$, $B=B(h(t),s)$, $C=C(h(t),s)$, $D=D(h(t),s)$
and $u_0=S_0, v_0=I_0, w_0=R_0$.

We denote $h^*=-\mu v'_0(h_0)$, and for $0<T\leq\frac {h_0}{8(1+h^*)}$, set
$$H_T=\Big\{h\in C^1[0,T]: \, h(0)=h_0,  \ h'(0)=h^*, \, ||h'-h^*||_{C([0, T])}\leq 1\Big\},$$
$$U_T=\Big\{u\in C([0, +\infty)\times[0,T]): \, u(s,0)=u_0(s),\,  \|u-u_0\|
_{L^\infty([0, +\infty)\times [0, T])}\leqslant 1 \Big\},$$
\begin{eqnarray*}
&V_T&=\Big\{v\in C([0, \infty)\times[0,T]): \,v(s,t)\equiv 0 \,
      \textrm{for}\, s\geq h_0, 0\leq t\leq T,\\
& &v(s,0)=v_0(s)\, \textrm{for}\, 0\leq s\leq h_0, \
\|v-v_0\|_{L^\infty([0, \infty)\times [0, T])}\leqslant 1\Big\},
\end{eqnarray*}
\begin{eqnarray*}
&W_T&=\Big\{w\in C([0, \infty)\times[0,T]): \,w(s,t)\equiv 0 \,
      \textrm{for}\, s\geq h_0, 0\leq t\leq T,\\
& &w(s,0)=w_0(s)\; \textrm{for}\; 0\leq s\leq h_0, \
\|w-w_0\|_{L^\infty([0, \infty)\times [0, T])}\leqslant 1\Big\}.
\end{eqnarray*}

Noticing the fact that for $h_1, h_2\in
H_T$, due to $h_1(0)=h_2(0)=h_0$, we have
\begin{equation}
\label{Bc}
 \|h_1-h_2\|_{C([0,T])}\leq  T||h'_1-h'_2||_{C([0, T])},
\end{equation}
it is not difficult to
see that $\Gamma_T:=U_T\times V_T\times W_T\times H_T$ is a complete metric space with
the metric
$$\mathcal{D}((u_1, v_1, w_1; h_1),(u_2, v_2, w_2; h_2))=\|u_1-u_2\|_{L^\infty ([0, +\infty)\times[0,T])}$$
$$+\|v_1-v_2\|_{L^\infty([0,\infty)\times [0,T])}+\|w_1-w_2\|_{L^\infty([0, \infty)\times [0,T])}+
\|h'_1-h'_2\|_{C([0,T])}.$$

Next, we shall prove the existence and uniqueness result by using
the contraction mapping theorem. Applying standard $L^p$ theory and the Sobolev imbedding theorem
\cite{LSU}, we can find that for any $(u, v, w; h)\in\Gamma_T$, the
following initial boundary value problem
\begin{eqnarray}
\left\{
\begin{array}{lll}
\tilde u_{t}-Ad_1\Delta_s \tilde u-(Bd_1+h'C+Dd_1)\tilde u_s=b-\beta uv-\mu _1u,\; &
s>0,\; t>0,  \\[1mm]
\tilde v_{t}-Ad_2\Delta_s
\tilde v-(Bd_2+h'C+Dd_2)\tilde v_s=\beta uv-\mu _2v-\alpha v,\; &
0<s<h_0,\; t>0, \\[1mm]
\tilde w_{t}-Ad_3 \Delta_s \tilde w-(Bd_3+h'C+Dd_3)\tilde w_s=\alpha v-\mu _3w,\; &
0<s<h_0,\; t>0, \\[1mm]
\tilde u_s(0, t)=\tilde v_s(0, t)=\tilde w_s(0,t)=0, &t>0,\\[1mm]
\tilde v(s,t)=\tilde w(s,t)=0,&s\geq h_0,\, t>0,\\[1mm]
\tilde u(s,0)=u_{0}(s),\, \tilde v(s,0)=v_0(s),\, \tilde w(s,0)=w_0(s),&s\leq 0
\end{array} \right.
\label{Bd}
\end{eqnarray}
admits a unique solution $$(\tilde{u}, \tilde v, \tilde w )\in [C^{1+\gamma, (1+\gamma)/2
}([0, +\infty)\times [0, T])]^3$$
 and
\begin{eqnarray}
&&\|\tilde{u}\|_{C^{1+\gamma,(1+\gamma)/2}([0,+\infty)\times[0,T])}\leqslant K_1,\label{Be}\\[1mm]
&&\|\tilde{v}\|_{C^{1+\gamma,(1+\gamma)/2}([0,h_0]\times[0,T])}\leqslant  K_1,\label{Be1}\\[1mm]
&&\|\tilde{w}\|_{C^{1+\gamma,(1+\gamma)/2}([0,h_0]\times[0,T])}\leqslant K_1,\label{Be2}
\end{eqnarray}
where $K_1$ is a constant depending on
$\gamma, h_0$, $\|S_0\|_{C^{2}[0, +\infty)}, \|I_0\|_{C^{2}[0, h_0]}$ and $\|R_0\|_{C^{2}[0, h_0]}.$

Now, we define $\tilde{h}(t)$ by the sixth equation in
(\ref{Bb}) as the following:
\begin{equation}
\label{Bf}
\tilde{h}(t)=h_0-\mu\int^t_0\tilde{v}_s(h_0, \tau)\textrm{d}\tau,
\end{equation}
then we have $\tilde{h}'(t)=-\mu \tilde{v}_s(h_0,t)$,
$\tilde{h}(0)=h_0$ and $\tilde{h}'(0)=-\mu
v'_0(h_0)=h^*$. Hence $\tilde{h}'(t)\in C^{\gamma/2}([0,T])$ with
\begin{equation}
\label{Bg}
\|\tilde{h}'(t)\|_{C^{\gamma/2}([0,T])}\leq K_2:=\mu K_1.
\end{equation}

In what follows, we define a map
$$\mathcal{F}:\ \Gamma_{T}\longrightarrow [C([0,+\infty)\times[0,T])]^3\times
C^1([0,T])$$  by
$\mathcal{F}(u(s,t), v(s,t), w(s,t); h(t)) = (\tilde u(s,t), \tilde v(s,t), \tilde w(s,t); \tilde{h}(t))$.
It is obvious that $(u(s,t), v(s,t), w(s,t); h(t))\in \Gamma_T$ is a fixed point of
$\mathcal{F}$ if and only if it solves (\ref{Bb}).

Similarly as in \cite{DL}, there is a $T>0$ such that $\mathcal{F}$ is a contraction mapping in $\Gamma_{T}$. It follows from
 the contraction mapping theorem that
there is a $(u(s,t), v(s,t), w(s,t); h(t))$ in $\Gamma_T$ such that
$$\mathcal{F}(u(s,t), v(s,t), w(s,t); h(t))=(u(s,t), v(s,t), w(s,t); h(t)).$$
 In other words, $(u(s,t), v(s,t), w(s,t); h(t))$ is the solution of the
problem (\ref{Bb}) and thereby $(S(r,t), I(r,t), R(r,t); h(t))$ is the
solution of the problem (\ref{a3}). Moreover, by using the Schauder
estimates, we have additional regularity of the solution, $h(t)\in
C^{1+\gamma/2}([0,T])$,  $S\in
C^{2+\gamma,1+\gamma/2}((0, +\infty )\times (0,T])$ and $I, R\in
C^{2+\gamma,1+\gamma/2}((0, h(t))\times (0,T])$. Thus $(S(r, t), I(r, t), R(r, t); h(t))$
is the classical solution of the problem (\ref{a3}). \epf

To show that the local solution obtained in Theorem 2.1 can be
extended to all $t>0$, we need the following estimate.

\begin{lem} Let $(S, I, R; h)$ be a bounded solution to problem \eqref{a3} defined for $t\in (0,T_0)$ for some $T_0\in (0, +\infty]$.
Then there exist constants $C_1$ and $C_2$ independent of $T_0$ such
that
\[
0<S(r, t)\leq C_1\; \mbox{ for } 0\leq r<
+\infty,\; t\in (0, T_0). \]
\[
0<I(r, t), R(r,t)\leq C_2\; \mbox{ for } 0\leq r<
h(t),\; t\in (0,T_0). \]
\end{lem}
\bpf
It is easy to see that $S\geq 0, I\geq 0$  and $R\geq 0$ in $[0, +\infty)\times
[0, T_0)$ as long as the solution exists.

Using the strong maximum principle to the equations in $[0, h(t)]\times
[0, T_0)$, we immediately obtain
 $$S(r, t), I(r, t), R(r,t)>0\;\; \textrm{for} \ 0\leq r<h(t),\,0< t<T_0.$$

 The upper bounds of the solution are followed from the maximum principle, we omit the proof here.
\epf

The next lemma shows that the free boundary for problem (\ref{a3}) is strictly monotone increasing.

\begin{lem}
Let $(S, I, R; h)$ be a solution to problem \eqref{a3} defined for
$t\in (0,T_0)$ for some $T_0\in (0, +\infty]$.
Then there exists a constant $C_3$  independent of $T_0$ such that
\[
 0<h'(t)\leq C_3 \; \mbox{ for } \; t\in (0,T_0). \]
\end{lem}
\bpf Using the Hopf Lemma to the equation of $I$ yields that
 $$I_r(h(t), t)<0 \ \;\; \textrm{for} \ 0< t<T_0.$$
 Hence $h'(t)>0$ for $t\in (0,T_0)$ from the Stefan condition.

Next we show that $h'(t)\leq C_3$ for all $t\in (0,T_0)$ and
some $C_3$ independent of $T_0$. As in \cite{LIN} , we define
$$\Omega=\Omega_M: =\{(r, t):h(t)-M^{-1}<r<h(t),\, \, 0<t<T_0\}$$ and
construct an auxiliary function
$$w(r, t):=C_2[2M(h(t)-r)-M^2(h(t)-r)^2].$$
We will choose $M$ so that $w(r, t)\geq I(r, t)$ holds over $\Omega$.

Direct calculations show that, for $(r, t)\in\Omega$,
$$w_t=2C_2Mh'(t)(1-M(h(t)-r))\geq 0,$$
$$-\Delta w=2C_2M^2,\quad \beta SI-(\mu_2+\alpha )I\leq \beta C_1C_2,$$
and then $$w_t-d_2 \Delta w\geq 2 d_2 C_2M^2\geq \beta C_1C_2\mbox{  in }
\Omega $$ if $M^2\geq \frac {\beta C_1}{2d_2}$. On the other hand, we have
$$w(h(t)-M^{-1},t)=C_2\geq I(h(t)-M^{-1}, t), \quad
w(h(t), t)=0=I(h(t), t).$$
Hence, if we can choose $M$ such that
$I_0(r)\leq w(r, 0)$ for $r\in [h_0-M^{-1}, h_0]$, then we can apply
the maximum principle to $w-I$ over $\Omega$ to deduce that
$I(r, t)\leq w(r, t)$ for $(r, t)\in\Omega$. It would then follow that
\[I_r( h(t), t)\geq w_r(h(t), t)=-2MC_2,\; h'(t)=-\mu I_r(
h(t), t)\leq C_3:=2MC_2\mu.\]

To complete the proof, we only have to find some $M$ independent of
$T_0$ such that $I_0(r)\leq w(r, 0)$ for $r\in [h_0-M^{-1}, h_0]$. We
calculate
\[
w_r(r, 0)=-2C_2M[1-M(h_0-r)]\leq -C_2M \mbox{ for } r\in
[h_0-(2M)^{-1}, h_0].
\]
Then upon choosing
\[
M:=\max\left\{ \sqrt{\frac{\beta C_1}{2d_2}},
\frac{4\|I_0\|_{C^1([0,h_0])}}{3C_2}\right\},
\]
we have
\[ w_r(r, 0)\leq -MC_2\leq -\frac 43||I_0||_{C^1}\leq I_0'(r) \mbox{ for } r\in
[h_0-(2M)^{-1}, h_0].
\]
Since $w(h_0, 0)=I_0(h_0)=0$, the above inequality implies that
\[ w(r, 0)\geq I_0(r) \mbox{ for } r\in
[h_0-(2M)^{-1}, h_0].
\]
Moreover, for $r\in [h_0-M^{-1}, h_0-(2M)^{-1}]$, we have
\[w(r, 0)\geq \frac 3 4 C_2,\; I_0(r)\leq
\|I_0\|_{C^1([0, h_0])}M^{-1}\leq \frac 34 C_2.
\]
Therefore $I_0(r)\leq w(r, 0)$ for $r\in [h_0-M^{-1}, h_0]$.
 This completes the proof.
 \epf

\begin{thm} The solution of the problem \eqref{a3} exists and is
unique for all $t\in (0,\infty)$.
\end{thm}
\bpf
It follows from the uniqueness of the solution that there is a number $T_{max}$
such that $[0, T_{max})$ is the maximal time interval in which the
solution exists. Now we prove that $T_{max}=\infty$ by the contradiction argument. Assume that
$T_{max}<\infty$. Then it follows from Lemma 2.2 that there exist $C_1, C_2$ and $C_3$
independent of $T_{max}$ such that for $t\in [0, T_{\max})$ and
$r\in [0, h(t)]$,
\[ 0\leq S(r,t)\leq C_1, \;(r,t)\in [0, +\infty)\times [0, T_{\max}),\]
\[ 0\leq I(r,t),\, R(r, t)\leq C_2, \;(r,t)\in [0, h(t)]\times [0, T_{\max}),\]
\[ h_0\leq h(t)\leq h_0+C_3 t, \ 0\leq h'(t)\leq C_3,\; t\in [0, T_{\max}).\]
We now fix $\delta_0\in (0, T_{max})$ and $M>T_{max}$. Then by the standard
parabolic regularity, we can find $C_4>0$ depending only on
$\delta_0$, $M$, $C_1$, $C_2$ and $C_3$ such that
 $$||S(\cdot, t)||_{C^{1+\gamma}[0, +\infty)},\ ||I(\cdot, t)||_{C^{1+\gamma}[0, h(t)]},\
 ||R(\cdot, t)||_{C^{1+\gamma}[0, h(t)]}\leq
C_4$$
 for $t\in [\delta_0, T_{\max})$.  It then follows from the
proof of Theorem 2.1 that there exists a $\tau>0$ depending only on
$C_i(i=1, 2, 3, 4)$ such that the solution of  problem \eqref{a3}
with initial time $T_{max}-\tau/2$ can be extended uniquely to the
time $T_{max}-\tau /2+\tau$. But this contradicts the assumption and thereby the proof is complete. \epf

\begin{rmk}
It follows from the uniqueness of the solution to \eqref{a3} and
some standard compactness arguments that the unique solution $(S,I,R,h)$
depends continuously on the parameters appearing in \eqref{a3}. This
fact will be used in the following sections hereafter.
\end{rmk}

We next decide when the transmission of diseases is spreading or vanishing. We
need to divide our discussion into two cases: $R_0<1$ and $R_0>1$.

\section{The case $ R_{0}<1 $}
It follows from Lemma 2.3 that $r=h(t)$ is monotonic increasing and therefore
there exists  $h_\infty\in (0, +\infty]$ such that $\lim_{t\to +\infty} \ h(t)=h_\infty$.
The following theorem shows that the transmission of diseases is vanishing in the case that $R_0<1$.
\begin{thm} If $R_0(=:\frac{b \beta}{\mu_1(\mu_2+\alpha)})<1$, then  $\lim_{t\to
+\infty} \ ||I(\cdot, t)||_{C([0, h(t)])}=0$ and $h_\infty<\infty$. Moreover, $\lim_{t\to
+\infty} \ ||R(\cdot, t)||_{C([0, h(t)])}=0$ and $\lim_{t\to +\infty} \ S(r, t)=\frac
{b}{\mu_1}$ uniformly in any bounded subset of $[0, \infty)$.
\end{thm}
\bpf
It follows from the comparison principle that $S(r, t)\leq \overline S(t)$
for $r\geq 0$ and $t\in(0, +\infty)$,
where
$$\overline S(t) :=\frac{b}{\mu_1}+(||S_0||_\infty-\frac{b}{\mu_1}) e^{-\mu_1 t},$$
which is the solution of the problem
\begin{equation}\label{ode}
\frac {d \overline S}{d t}= b-\mu_1 \overline S,\quad t>0;\;  \overline S(0)=||S_0||_\infty.
\end{equation}
Since $\lim_{t\to\infty}\overline
S(t)= \frac b\mu_1$, we deduce that
\[
\mbox{ $\limsup_{t\to +\infty} \ S(r, t)\leq \frac b\mu_1$ uniformly for
$r\in [0,\infty)$.}
\]

Recalling the condition $R_0<1$, there exists $T_0$ such $S(r, t)\leq \frac b\mu_1\frac {1+R_0}{2R_0}$
in $[0, \infty)\times [T_0, +\infty)$.
Now $I(r,t)$ satisfies \begin{eqnarray}
\left\{
\begin{array}{lll}
I_t-d_2 \Delta I\leq [\frac
{\beta b}{\mu_1}\frac {1+R_0}{2R_0}-\mu _2-\alpha ]I(r, t),\; &
0<r<h(t),\; t>T_0, \\[1mm]
I(r,t)=0,&r=h(t),\, t>0,\\[1mm]
I(r,T_0)>0,\,  &0\leq r\leq h(T_0).
\end{array} \right.
\label{am13}
\end{eqnarray}
Therefore $||I(\cdot, t)||_{C([0, h(t)])}\to 0$
as $t\to \infty$, since that $\beta \frac b\mu_1\frac {1+R_0}{2R_0}-\mu _2-\alpha <0$.
We then have $||R(\cdot, t)||_{C([0, h(t)])}\to 0$
as $t\to \infty$ from the third equation of (\ref{a3}).

 Next we show that $h_\infty <+\infty$. In fact, direct calculation yields
\begin{eqnarray*}& &\frac{\textrm{d}}{\textrm{d} t}\int_0^{h(t)}r^{n-1}I(r, t)\textrm{d}r\\[1mm]
&=&\int_0^{h(t)}r^{n-1}I_t(r,t)\textrm{d}r+h'(t)h^{n-1}(t)I(h(t), t)\\[1mm]
&=&\int_0^{h(t)}d_2 r^{n-1}\Delta I\textrm{d}r+\int_0^{h(t)}I(r, t)(\beta S(r, t)-\mu _2-\alpha)r^{n-1}\textrm{d}r\\[1mm]
&=&\int_0^{h(t)}d_2(r^{n-1}I_r(r,t))_r\textrm{d}r+\int_0^{h(t)}I(r, t)(\beta S(r, t)-\mu _2-\alpha)r^{n-1}\textrm{d}r\\[1mm]
&=&-\frac{d_2}{\mu}h^{n-1}h'(t)+\int_0^{h(t)}I(r, t)(\beta S(r, t)-\mu _2-\alpha)r^{n-1}\textrm{d}r.
\end{eqnarray*}

Integrating from $T_0$ to $t\,(>T_0)$ yields
\begin{eqnarray}
& &\int_0^{h(t)}r^{n-1}I(r, t)\textrm{d}r = \int ^{h(T_0)}_0 r^{n-1}I(r, T_0)dr
+\frac {d_2}{n\mu}h^n(T_0)-\frac {d_2}{n\mu}h^n(t)\nonumber \\[1mm]
& &\quad +\int_{T_0}^t\int_0^{h(s)}I(r, s)(\beta S(r, s)-\mu _2-\alpha)r^{n-1}drds,
\quad t\geq T_0.\label{k1}
\end{eqnarray}
Since $0<S(r, t)\leq\frac b\mu_1\frac {1+R_0}{2R_0}$ for $r\in [0,h(t))$ and $t\geq T_0$, we have
$$\beta S(r, t)-\mu _2-\alpha\leq 0 \mbox{ for $t\geq T_0$},$$
$$\int_0^{h(t)}r^{n-1}I(r, t)\textrm{d}r \leq \int ^{h(T_0)}_0 r^{n-1}I(r, T_0)dr +\frac
{d_2}{n\mu}h^n(T_0)-\frac {d_2}{n\mu}h^n(t) \mbox{ for $t\geq T_0$},$$
which in turn gives that $h_\infty<\infty$.

Then it follows from the first equation of (\ref{a3}) that $\lim_{t\to +\infty} \ S(r, t)=\frac
{b}{\mu_1}$ uniformly in any bounded subset of $[0, \infty)$.
\epf

\section{The case $ R_{0}>1 $}
In order to study the case that the reproduction number $R_0>1$, and for later applications, we
 need a comparison principle, which can be used to estimate $S(r,t)$, $I(r,t)$, $R(r,t)$ and the free
boundary $r=h(t)$. As in \cite{DL}, the following comparison lemma can be obtained analogously.

\begin{lem} Suppose that $T\in (0,\infty)$, $\overline h\in
C^1([0,T])$, $\overline S\in C([0, \infty)\times [0,T])\cap C^{2,1}((0, \infty)\times (0,T])$,
$\overline I, \overline R\in C(\overline D_T^*)\cap C^{2,1}(D_T^*)$
with $D_T^*=\{(r,t)\in\R^2: 0<r<\overline h(t),0<t\leq T\}$, and
\begin{eqnarray*}
\left\{
\begin{array}{lll}
\bar{S}_t-d_1\Delta\bar{S}\geq b-\mu_{1}\bar{S},\; &0<r,\ 0<t\leq T, \\[1mm]
\bar{I}_t-d_2\Delta \bar{I}\geq(\beta\bar{S}-\mu_{2}
-\alpha)\bar{I},\; &0<r<\bar{h}(t),\ 0<t\leq T, \\[1mm]
\bar{R}_t-d_3\Delta \bar{R}\geq\alpha\bar{I}-\mu_{3}
\bar{R},\; &0<r<\bar{h}(t),\ 0<t\leq T, \\[1mm]
\bar{S}_{r}(0,t)\geq 0,\ \bar{I}_{r}(0,t)\geq0,\ \bar{R}_{r}(0,t)\geq0,\; & 0<t\leq T, \\[1mm]
\bar{I}(r,t)=\bar{R}(r,t)=0,\; &r\geq\bar{h}(t),\ 0<t\leq T, \\[1mm]
\bar h'(t)\geq -\mu \bar{I}_{r}(\bar h(t),t),\ \bar h(0)>h_{0},\; &0<t\leq T,\\[1mm]
\bar{S}(r,0)\geq S_{0}(r),\ \bar{I}(r,0)\geq I_{0}(r),\ \bar{R}(r,0)\geq R_{0}(r),\; &0\leq r\leq h_0.
\end{array} \right.
\end{eqnarray*}
Then the solution $(S,I,R;h)$ of free boundary problem \eqref{a3}
satisfies
\[\mbox{ $S(r,t)\leq \overline S(r,t)$, $h(t)\leq\overline h(t)$ for $r\in (0, \infty)$ and $t\in (0, T]$,}\]
\[\mbox{ $I(r,t)\leq \overline I(r,t)$, $R(r,t)\leq \overline R(r,t)$ for $r\in (0,h(t))$ and $t\in (0, T]$.}\]
\end{lem}

Next we show that if $h_0$ and $\mu$ are sufficiently small, the disease is vanishing for the case $R_0>1$.
\begin{thm}
If $R_0(:=\frac{b \beta}{\mu_1(\mu_2+\alpha)})>1$,
$h_0\leq \min \left\{ \sqrt{\frac{d_2}{16k_{0}}},\ \sqrt{\frac{d_2}{16\alpha}} \right\}$
and $\mu\leq \frac{d}{8M}$, then $h_\infty<\infty$.
Where $k_0=\beta C_{1}-\mu_{2}-\alpha>0$, $C_{1}=\max\left\{||S_{0}||_\infty,\ \frac{b}{\mu_{1}}
\right\}$, and $M=\frac{4}{3}\max\left\{||I_{0}||_\infty,\ ||R_{0}||_\infty \right\}$.
\end{thm}
\bpf
We are going to construct a suitable upper solution to
\eqref{a3} and then apply Lemma 4.1. As in \cite{DL}, we define
$$\bar{S}(r,t)=C_{1},\qquad \qquad \qquad\qquad\qquad\qquad\qquad$$
\begin{eqnarray*}
\bar{I}=
\left\{
\begin{array}{lll}
M e^{-\gamma t}V(r/\overline h (t)),\; & 0\leq r\leq \overline h(t),\\[1mm]
0,\; &r> \overline h(t),
\end{array} \right.
\end{eqnarray*}
\begin{eqnarray*}
\bar{R}=
\left\{
\begin{array}{lll}
M e^{-\gamma t}V(r/\overline h (t)),\; & 0\leq r\leq \overline h(t), \\[1mm]
0,\; &r> \overline h(t),
\end{array} \right.
\end{eqnarray*}
and $$\overline h (t)=2h_0(2- e^{-\gamma t}), \  t\geq 0;
 \quad V(y)=1-y^{2}, \ 0\leq y\leq 1,$$
where $C_{1}=\max\left\{||S_{0}||_\infty,\ \frac{b}{\mu_{1}}
\right\}$, $\gamma$ and $M$ are
positive constants to be chosen later.

Denoting $k_0=\beta C_{1}-\mu_{2}-\alpha$, we have $k_0>0$ since $R_0>1$. Direct computations yield
\begin{eqnarray*}
& &\bar{S}_t-d_1\Delta\bar{S}=0\geq b-\mu_{1}\bar{S},\\[1mm]
& &\bar{I}_t-d_2\Delta \bar{I}-(\beta\bar{S}-\mu_{2}-\alpha)\bar{I}\\[1mm]
& &=\bar{I}_t-d_2\Delta \bar{I}-k_{0}\bar{I}\\[1mm]
& &=M e^{-\gamma t}[-\gamma V-r\overline h '\overline h^{-2}V'-d_{2}
\overline h^{-2}V''-d_2\frac{n-1}{r}\overline h^{-1}V'-k_{0}V]\\[1mm]
& &\geq M e^{-\gamma t}[\frac{d_{2}}{8h_0^2}-\gamma-k_{0}],\\[1mm]
& &\bar{R}_t-d_3\Delta \bar{R}-(\alpha\bar{I}-\mu_{3}
\bar{R})\geq M e^{-\gamma t}[\frac{d_3}{8h_0^2}-\gamma-\alpha]
\end{eqnarray*}
for all $0<r<\overline h(t)$ and $t>0$. On the other hand, we have
$\overline h'(t)=2h_0 \gamma e^{-\gamma t}$ and
$-\mu \bar{I}_r(\overline h(t), t)=2M \mu \overline h^{-1}(t)e^{-\gamma t}$.
Moreover, it follows that $\bar{S}(r, 0)\geq S_{0}(r)$, $\bar{I}(r, 0)=M(1-\frac{r^{2}}{4h_{0}^{2}})\geq\frac{3}{4}M$,
 $\bar{R}(r, 0)=M(1-\frac{r^{2}}{4h_{0}^{2}})\geq\frac{3}{4}M$ for
$r\in [0, h_0]$. Noting that $\overline h(t)\leq4h_0$, we now choose
$$M=\frac{4}{3}\max\left\{||I_{0}||_\infty,\ ||R_{0}||_\infty \right\}$$
and take
$$\gamma= \frac{d}{16h^2_{0}},\ \mu\leq \frac{d}{8M},$$
where $d:=\min\{d_1, d_2\}$ and $h_0\leq \min \left\{ \sqrt{\frac{d}{16k_{0}}},\
 \sqrt{\frac{d}{16\alpha}} \right\}$.
Then we have
\begin{eqnarray*}
\left\{
\begin{array}{lll}
\bar{S}_t-d_1\Delta\bar{S}\geq b-
\beta\bar{S}\underline{I}-\mu_{1}\bar{S},\; &0<r,\ t>0, \\[1mm]
\bar{I}_t-d_2\Delta \bar{I}\geq(\beta\bar{S}-\mu_{2}
-\alpha)\bar{I},\; &0<r<\overline h(t),\ t>0, \\[1mm]
\bar{R}_t-d_3\Delta \bar{R}\geq\alpha\bar{I}-\mu_{3}
\bar{R},\; &0<r<\overline h(t),\ t>0, \\[1mm]
\bar{S}_{r}(0,t)=0,\ \bar{I}_{r}(0,t)\geq0,\ \bar{R}_{r}(0,t)\geq0,\; & t>0, \\[1mm]
\bar{I}(r,t)=\bar{R}(r,t)=0,\; &r\geq\overline h(t),\ t>0, \\[1mm]
\overline h'(t)\geq -\mu \bar{I}_{r}(\overline h(t),t),\ \overline h(0)=2h_{0}>h_{0},\; &t>0,\\[1mm]
\bar{S}(r,0)\geq S_{0}(r),\ \bar{I}(r,0)\geq I_{0}(r),\ \bar{R}(r,0)\geq R_{0}(r),\; &r\geq 0.
\end{array} \right.
\end{eqnarray*}
Hence we can apply Lemma 4.1 to conclude that  $h(t)\leq\overline h(t)$
for $t>0$. Therefore, we have $h_\infty\leq \lim_{t\to\infty}
\overline h(t)=4h_0<\infty$.
 \epf

For the case that $R_0>1$, we next prove that if $h_0$ is suitably large,
the disease is spreading.
\begin{lem}  If $h_\infty<\infty$, then $\lim_{t\to
+\infty} \ ||I(\cdot, t)||_{C([0, h(t)])}=0$. Moreover, we have $\lim_{t\to
+\infty} \ ||R(\cdot, t)||_{C([0, h(t)])}=0$ and $\lim_{t\to +\infty} \ S(r,t)=\frac
{b}{\mu_1}$ uniformly in any bounded subset of $[0, \infty)$.
\end{lem}
\bpf
Assume  $\limsup_{t\to
+\infty} \ ||I(\cdot, t)||_{C([0, h(t)])}=\delta>0$ by contradiction.
Then there exists a sequence $(r_k, t_k )$ in $[0, h(t))\times(0,\infty)$
such that $I( r_k,t_k)\geq \delta /2$ for all $k \in \mathbb{N}$,
and $t_k\to \infty$ as $k\to \infty$.
Since that $0\leq r_k<h(t)<h_\infty<\infty$,
we then have that a subsequence of $\{r_n\}$ converges to $r_0\in [0, h_\infty)$.
Without loss of generality,
we assume $r_k\to r_0$ as $k\to \infty$.

Define $S_k(r,t)=S(r,t_k+t)$, $I_k(r,t)=I (r,t_k+t)$ and $R_k(r,t)=R(r,t_k+t)$ for
$r\in (0, h(t_k+t)),t\in (-t_k, \infty)$.
It follows from the parabolic regularity that  $\{(S_k, I_k, R_k)\}$ has a subsequence $\{(S_{k_i}, I_{k_i}, R_{k_i})\}$ such that
$(S_{k_i}, I_{k_i}, R_{k_i})\to (\tilde S, \tilde I, \tilde R)$ as $i\to \infty$ and $(\tilde S, \tilde I, \tilde R)$ satisfies
\begin{eqnarray*}
\left\{
\begin{array}{lll}
\tilde S_t-d_1 \Delta \tilde S=b-\beta \tilde S \tilde I-\mu_1\tilde S,\;
      & \ 0<r<h_\infty,\ t\in (-\infty, \infty), \\[1mm]
\tilde I_t-d_2 \Delta \tilde I=\beta \tilde S \tilde I-\mu_2\tilde I-\alpha \tilde I,\;
      &\ 0<r<h_\infty, \ t\in (-\infty, \infty),\\[1mm]
\tilde R_t-d_3 \Delta \tilde R=\alpha \tilde I-\mu_3 \tilde R,\;
&\ 0<r<h_\infty,\ t\in (-\infty, \infty).
\end{array} \right.
\end{eqnarray*}
Since $\tilde I(r_0,0 )\geq \delta/2$,
we have $\tilde I>0$ in $ [0, h_\infty)\times(-\infty, \infty)$.
Recalling that $(\beta \tilde S -\mu_2-\alpha)$ is bounded by
$M:=\beta \max\{\frac {b}{\mu_1}, ||S_0||_{L^\infty}\}+\mu_2+\alpha$.
Applying Hopf lemma to the equation
$\tilde I_t-d_2 \Delta \tilde I\geq -M\tilde I$ at the point $(0, h_\infty)$ yields
that $\tilde I_r(h_\infty, 0 )\leq -\sigma_0$ for some $\sigma_0>0$.

On the other hand, $h(t)$ is increasing and bounded.
Moreover, for any $0<\alpha <1$, there exists a constant $\tilde C$,
which depends on $\alpha, h_0, \|I_0\|_{C^{1+\alpha}[0, h_0]}$ and $h_\infty$, such that
\begin{eqnarray}
\|I\|_{C^{1+\alpha,(1+\alpha)/2}([0, h(t))\times[0, \infty))}
+\|h\|_{C^{1+\alpha/2}([0,\infty))}\leq\tilde C.\label{Bg1}
\end{eqnarray}
In fact, let us straighten the
free boundary in a way different from that in Theorem 2.1. Define
$$s=\frac{h_0r}{h(t)},\ u(s,t)=S(r,t), \ v(s,t)=I(r,t), \ w(s,t)=R(r,t),$$
then direct calculations yield that
$$I_t=v_t-\frac{h'(t)}{h(t)}sv_s,\ I_r=\frac{h_0}{h(t)}v_s,\
  \Delta_r I=\frac{h^2_0}{h^2(t)}\Delta_s v,$$
therefore, $v(s, t)$ satisfies
\begin{eqnarray}
\left\{
\begin{array}{lll}
v_{t}-d_2\frac{h^2_0}{h^2(t)}\Delta_s v-\frac{h'(t)}{h(t)}sv_s=
  v(\beta u -\mu_2-\alpha),\; &0<s<h_0,\; t>0, \\[1mm]
v_s(0,t)=v(h_0,t)=0,&t>0,\\[1mm]
v(s,0)=I_0(s)\geq 0,\; & 0\leq s\leq h_0.
\end{array} \right.
\label{Bb1}
\end{eqnarray}
This transformation changes the free boundary $r=h(t)$ to the fixed
line $s=h_0$, at the expense of making the equation more
complicated. It follows from Lemmas 2.2 and 2.3 that
$$||v(\beta u -\mu_2-\alpha)||_{L^\infty}\leq M_1,\ ||\frac{h'(t)}{h(t)}s||_{L^\infty}\leq M_3.$$
Applying standard $L^p$ theory  and then the Sobolev imbedding
theorem (\cite{LSU}), we obtain that
\begin{eqnarray*}
\|v\|_{C^{1+\alpha,(1+\alpha)/2}([0,h_0]\times[0,\infty))}\leq  C_4,
\end{eqnarray*}
where $C_4$ is a constant depending on
$\alpha, h_0, M_1, M_2, M_3$ and $\|I_0\|_{C^{2}[0, h_0]}$. This immediately leads to (\ref{Bg1}).

Since $\|h\|_{C^{1+\alpha/2}([0,\infty))}\leq\tilde C$ and $h(t)$ is bounded,
we then have $h'(t)\to 0$ as $t\to \infty$, that is,
$I_r(h(t_k),t_k)\to 0$ as $t_k\to \infty$ by the free boundary condition.
Moreover, it follows from the inequality $\|I\|_{C^{1+\alpha,(1+\alpha)/2}([0, h(t))
\times[0, \infty))}\leq \tilde C$ that
$I_r(h(t_k),t_k+0)=(I_k)_r(h(t_k),0)\to \tilde I_r(h_\infty,0)$ as $k\to \infty$,
which leads to a contradiction to the fact that  $\tilde I_r(h_\infty,0)\leq -\sigma_0<0$.
 Thus $\lim_{t\to +\infty} \ ||I(\cdot,t)||_{C([0,h(t)])}=0$, and thereby $\lim_{t\to
+\infty} \ ||R(\cdot, t)||_{C([0, h(t)])}=0$ and $\lim_{t\to +\infty} \ S(r,t)=\frac
{b}{\mu_1}$ uniformly in any bounded subset of $[0, \infty)$.
\epf

Let $\lambda_1(R)$ be the principal eigenvalue of the operator $-\Delta$ in $B_R$ (open ball with radius $R$) subject to
homogeneous Dirichlet boundary condition.
It is well-known that $\lambda_1(R)$ is a strictly decreasing continuous function and
$$\lim_{R\to 0^+} \lambda_1(R)=+\infty \ \textrm{and}\ \lim_{R\to +\infty} \lambda_1(R)=0.$$

\begin{thm} If $R_0(:=\frac{b \beta}{\mu_1(\mu_2+\alpha)})>1$, then $h_\infty=\infty$ provided
that $h_0>h^*_0$, where $\lambda_1(h^*_0)=\frac{(\mu_2+\alpha)}{d_2}(R_0-1)$.
\end{thm}
\bpf Assume that  $h_\infty<+\infty$ by contradiction. It follows from Lemma 4.3 that $\lim_{t\to
+\infty} \ ||I(\cdot, t)||_{C([0, h(t)])}=0$. Moreover, $\lim_{t\to +\infty} \ S(r,t)=\frac
{b}{\mu_1}$ uniformly in the bounded subset $B_{h_0}$.
Therefore, for
$\varepsilon>0$, there exists $T^*>0$ such that $S(r, t)\geq \frac {b}{\mu_1}-\varepsilon$ for
$t\geq T^*, r\in [0, h(t))$.
We then have that $I(r,t)$ satisfies
  \begin{eqnarray}
\left\{
\begin{array}{lll}
I_{t}-d_2 \Delta I\geq (\beta (\frac {b}{\mu_1}-\varepsilon)-\mu_2-\alpha)I,\;
      & 0<r<h_0,\ t>T^*,\\[1mm]
I_r(0,t)=0,\; I(h_0,t)\geq 0,\quad & t>T^*,\\[1mm]
I(r,T^*)>0, &0\leq r<h_0.
\end{array} \right.
\label{fgg1}
\end{eqnarray}
It is easy to see that $I(r,t)$ has a lower solution $\underline I(r, t)$ satisfying
\begin{eqnarray}
\left\{
\begin{array}{lll}
\underline I_{t}-d_2 \Delta \underline I=
(\beta(\frac {b}{\mu_1}-\varepsilon)-\mu_2-\alpha)\underline I,\; & 0<r<h_0,\ t>T^*,\\[1mm]
\underline I_r(0,t)=0,\; \underline I(h_0,t)=0,\quad & t>T^*,\\[1mm]
\underline I(r,T^*)=I(r,T^*), &0\leq r<h_0.
\end{array} \right.
\label{fgg2}
\end{eqnarray}
Since $h_0>h^*_0$, we can choose $\varepsilon$ sufficiently small such that
$(\beta (\frac {b}{\mu_1}-\varepsilon)-\mu_2-\alpha)>d_2 \lambda_1(h_0)$,
it follows from well-known result that $\underline I$ is unbounded in
$(0, h_0)\times [T^*, \infty)$, which leads to a contradiction that $\lim_{t\to
+\infty} \ ||I(\cdot, t)||_{C([0, h(t)])}=0$.
 \epf

\section{Discussion}

In this paper, we have considered the SIR epidemic model
describing the transmission of diseases and examined the dynamical behavior
of the population $(S, I, R)$ with spreading front $r=h(t)$ determined by \eqref{a3}.
We have obtained the asymptotic behavior results.

The basic reproduction
number $R_0$ ($=\frac{b \beta}{\mu_1(\mu_2+\alpha)}$) is important but not unique factor to determine
whether the disease dies out or remains endemic. It is shown that
if $R_0<1$, vanishing always happens or the disease dies out (Theorem 3.1).
If $R_0>1$, spreading happens provided that $h_0$ is sufficiently
large (Theorem 4.4) and vanishing is possible provided that $h_0$ is small (Theorem 4.2).

We feel it is reasonable to conclude that \eqref{a3} is promising alternatives to \eqref{a1}
and \eqref{a2} for the modeling of
disease spreading, and there is still some works to do for the model \eqref{a3}. The first one is
that, what is the asymptotic spreading speed when spreading happens? Since there is no other choice except spreading and vanishing,
the second one is that we want to know the necessary condition for the disease to spread or to vanish.

\end{document}